\documentclass[11pt]{article}
\usepackage{graphicx}
\usepackage{amssymb}
\usepackage{amsmath}
\usepackage{amsfonts}
\usepackage{amsthm}

\voffset = - 1.2cm
\theoremstyle{plain}
\newtheorem{theorem}{Theorem}
\newtheorem{corollary}[theorem]{Corollary}
\newtheorem{proposition}[theorem]{Proposition}
\newtheorem{lemma}[theorem]{Lemma}
\newtheorem{example}{Example}

\theoremstyle{definition}
\newtheorem{definition}[theorem]{Definition}
\newtheorem{note}[theorem]{Note}

\numberwithin{equation}{section}
\numberwithin{theorem}{section}
\numberwithin{example}{section}

  2 
scaled\magstep3

\advance\topmargin by-1cm

\begin{document}

\centerline{\Large {\bf Mixed superposition rules and the Riccati hierarchy}} \vskip 0.75cm

\centerline{Janusz Grabowski and Javier de Lucas} \vskip 0.5cm

\centerline{Institute of Mathematics, Polish Academy of Sciences,}
\medskip
\centerline{\'Sniadeckich 8, P.O. Box 21, 00-956 Warszawa, Poland}

\vskip 1cm

\begin{abstract}
Mixed superposition rules, i.e., functions describing the general solution of a system of first-order
differential equations in terms of a generic family of particular solutions of first-order systems and some constants, are studied.
The main achievement is a generalization of the celebrated Lie-Scheffers Theorem, characterizing systems admitting a mixed superposition rule. This somehow unexpected result says that such systems are exactly Lie systems, i.e., they admit a standard superposition rule. This provides a new and powerful tool for finding Lie systems, which is
applied here to studying the Riccati hierarchy and to retrieving some known results in a more efficient and simpler way.
\end{abstract}

\section{Introduction}
{\it Lie systems} are systems of first-order differential equations whose general solutions can
be described in terms of generic families of particular solutions and some constants by a particular type
of maps: the {\it superposition rules} \cite{LS}-\cite{CarRamcinc}. The importance of these systems is due,
for instance, to their appearance in relevant physical and mathematical problems (see \cite{Dissertations} and
references therein). This has motivated a number of works devoted to the analysis of their properties and
applications \cite{FLV10}-\cite{CL08}.

The Lie-Scheffers Theorem \cite{LS,CGM07}, characterizing systems admitting a superposition rule, shows that
this property has geometrical roots and is rather exceptional. Although the superposition rules cannot be explicitly derived or effectively applied in many cases \cite{CL10SecOrd2}, their importance has motivated the search for generalizations \cite{CL10SecOrd2}-\cite{In65}.
For example, Vessiot pioneered the study of superposition rules for
second-order differential equations \cite{Ve95}, which was followed by Winternitz \cite{WintSecond} and by
Cari\~nena and coworkers \cite{CLR08,CL08Diss,CL10SecOrd2,CGL11}. Inselberg, in turn, proposed a concept of
superposition rules for operators \cite{In65}, while Winternitz and Schnider suggested a notion of
superposition rules for equations on supermanifolds \cite{HWA83,BecGagHusWin90}.

Our main aim in this note is the study of {\it mixed superposition rules}, whose original idea  was briefly
mentioned in \cite{CGM07,Dissertations} and implicitly used in
\cite{Dissertations,CLR08,Be07,Mil30,GGG11,GL99} in analyzing a number of important differential equations that appear in the physics and mathematical literature. On
the other hand, geometric properties of mixed superposition rules have not been thoroughly analyzed yet and we
hope that this work will fill in this gap.

As a main result, we characterize systems possessing a mixed superposition rule. More specifically, we prove
that a system of first-order differential equations admits a mixed superposition rule if and only if it is a
Lie system. This somehow unexpected theorem, called hereafter the {\it extended Lie-Scheffers Theorem}, is a
generalization of the Lie-Scheffers Theorem. In addition, it provides also a generalization of the {\it Lie's
condition} \cite{Dissertations}.

The extended Lie-Scheffers Theorem furnishes a new method to determine whether a system is a Lie system or
not: the search for a mixed superposition rule. As the latter is frequently easier than finding a standard
superposition rule, we easily recover several recent achievements about Riccati equations
\cite{PW,Dissertations}, Milne-Pinney equations
\cite{CLR08,WintSecond}, second-order Riccati equations \cite{CL10SecOrd2}, and Kummer-Schwarz equations \cite{CGL11} in a unified and simple way. Moreover, our approach can
be applied in other cases, e.g. the linearization of certain differential equations \cite{Be07}, where Lie systems
could also be employed.

Although we prove that the mixed superposition rules, like the standard superposition rules, can only appear in the context of Lie systems, the mixed superposition rules are much more versatile than the standard ones and allow us to express the general solution of a Lie system in a much broader variety of forms. Further, we show that mixed superposition rules can be used to investigate simultaneously different Lie systems.

Subsequently, in an analysis of mixed superposition rules, we demonstrate that there are some relations between properties of the involved Lie systems, which provides us with tools for the investigation of main features of one Lie system in terms of the characteristics of the others.

Our theoretical achievements provide a better understanding of geometrical properties of various particular equations of interest and furnish a framework for the study of questions of physical and mathematical relevance.
Special attention we paid to the analysis of the {\it Riccati hierarchy} \cite{Be07,CL08}, which is described
through an infinite family of Lie systems admitting mixed superposition rules. This permits us to combine the
general theory of Lie systems with our new methods in studying the whole Riccati hierarchy and other related problems \cite{PW,Dissertations,CLR08,CL08,CL10SecOrd2,Be07,Mil30,GL99}.

The structure of the paper goes as follows. Section 2 concerns the description of the basic notions to be used
throughout the paper. In Section 3 we present the motivation for working with mixed superposition rules and in
Section 4 we prove that mixed superposition rules can be understood as a certain type of flats connections that completes in depth an idea pointed out in \cite{CGM07}. Next, we provide in Section 5 a characterization of systems
admitting a mixed superposition rule and related results. In Section 6 we investigate a relevant type of mixed
superposition rules appearing in the literature. In Section 7 we use our theoretical results to
study the members of the Riccati hierarchy and other equations of physical interest. Finally, our main results and perspectives of a future work are summarized in Section 8.

\section{Preliminaries}\label{TDS}
For simplicity, we hereafter restrict ourselves to systems of differential equations on vector spaces and
assume geometric objects and functions to be real, smooth, and globally defined. This allows us to highlight
the main aspects of our work without discussing minor technical details.

Every system of first-order differential equations on $\mathbb{R}^{n_0}$,
\begin{equation}\label{LieSystem}
\frac{dx^i}{dt}=X^i(t,x),\qquad i=1,\ldots,n_0,
\end{equation}
is determined by the unique $t$-dependent vector field on $\mathbb{R}^{n_0}$,
\begin{equation}\label{VF}
X(t,x)=\sum_{i=1}^{n_0}X^i(t,x)\frac{\partial}{\partial x^i},
\end{equation}
whose integral curves (see \cite{CGL11} for details on this notion) are, up to a reparametrization, of the form $\gamma(t)=(t,x(t))$, with $x(t)$ being a
solution of (\ref{LieSystem}). Conversely, the $t$-dependent vector field $X$ determines its so-called {\it
associated system}, i.e. the system (\ref{LieSystem}) whose solutions describe its integral curves of the form
$\gamma:t\in\mathbb{R}\mapsto (t,x(t))\in\mathbb{R}\times\mathbb{R}^{n_0}$. This justifies the use of the
symbol $X$ for both, a $t$-dependent vector field and its associated system. Additionally, it can be proved that every $t$-dependent vector field $X$ is equivalent to a $t$-parametrized family $\{X^t\}_{t\in\mathbb{R}}$ of vector fields $X^t:\mathbb{R}^{n_0}\ni x\mapsto X(t,x)\in{\rm T}\mathbb{R}^{n_0}$  \cite{CGL11}.

\begin{definition}\label{LieSpan} Given a (finite or infinite) family $\mathcal{A}$ of vector fields on $\mathbb{R}^{n_0}$,
we denote with ${\rm Lie}(\mathcal{A})$ the smallest Lie algebra $V$ of vector fields on $\mathbb{R}^{n_0}$
containing $\mathcal{A}$. The  {\it minimal Lie algebra} $V^X$ of a $t$-dependent vector field $X$ is
$V^X={\rm Lie}(\{X^t\}_{t\in\mathbb{R}})$.
\end{definition}

Our work is mainly aimed to analyze mixed superposition rules. This concept represents a generalization of the
concept of a superposition rule.

\begin{definition} A {\it superposition rule} for a system $X$ is a function $\Phi:\mathbb{R}^{mn_0} \times \mathbb{R}^{n_0} \rightarrow\mathbb{R}^{n_0}$ of the form
\begin{equation}\label{super}
x=\Phi(x_{(1)},\ldots,x_{(m)};k_1,\ldots,k_{n_0})
\end{equation}
that allows us to write its general solution $x(t)$ as
\begin{equation}\label{superposition}
x(t)=\Phi(x_{(1)}(t),\ldots,x_{(m)}(t);k_1,\ldots,k_{n_0}),
\end{equation}
with $x_{(1)}(t),\ldots,x_{(m)}(t)$ being a generic family of particular solutions and $k_1,\ldots,k_{n_0}$
being constants related to the initial conditions for $x(t)$.
\end{definition}

The characterization of systems admitting a superposition rule is given by the celebrated {\it Lie-Scheffers
Theorem} \cite{LS,Dissertations}.

\begin{theorem} A system $X$ admits a superposition rule if and only if can be written in the form
\begin{equation}\label{LieDecom}
X^t=\sum_{\alpha=1}^rb_\alpha(t)Y_\alpha,
\end{equation}
for a family $Y_1,\ldots,Y_r$ of vector fields closing on an $r$-dimensional real Lie algebra and $t$-dependent functions $b_1(t),\ldots,b_r(t)$. In other words,
$X$ admits a superposition rule if and only if $V^X$ is finite-dimensional.
\end{theorem}

The finite-dimensional real Lie algebras of vector fields which are related to Lie systems are usually called
{\it Vessiot-Guldberg Lie algebras} in the literature \cite{Dissertations}. Apart from the Lie-Scheffers
Theorem, Lie also proved that if $X$ admits a superposition rule depending on $m$ particular solutions, then
$X$ possesses a Vessiot-Guldberg Lie algebra $V$ of dimension at most $m\cdot  n_0$. This is the referred to as
{\it Lie's condition} \cite{Dissertations}.

The Lie-Scheffers Theorem plays a central r\^ole in our work. Its geometric proof makes use of {\it diagonal prolongations}
\cite{CGM07,Dissertations}.

\begin{definition}\label{defProl} Given a $t$-dependent vector field $X$,
the $t$-dependent vector field $\widetilde X$ on $\mathbb{R}^{n_0(m+1)}$ of the form
\begin{equation}\label{prolongation}
\widetilde X=\sum_{a=0}^m\sum_{i=1}^{n_0}X^i(t,x_{(a)})\frac{\partial}{\partial x^i_{(a)}},
\end{equation}
is called the {\it diagonal prolongation} of $X$ to $\mathbb{R}^{n_0(m+1)}$.
\end{definition}

It is important to note that, given two $t$-independent vector fields $X_1,X_2$, we have $\widetilde{[X_1,X_2]}=[\widetilde X_1,\widetilde X_2]$. Hence, if $V$ is a Lie algebra of vector fields on
$\mathbb{R}^{n_0}$, the diagonal prolongations of its elements to $\mathbb{R}^{n_0m}$ span a Lie algebra
of vector fields isomorphic to $V$. Given a Lie algebra of vector fields $V$ on $\mathbb{R}^{n_0}$, we will
denote with $V_m$ the Lie algebra spanned by the diagonal prolongations of the elements of $V$ to
$\mathbb{R}^{n_0m}$.

\section{On the definition of mixed superposition rules}

It is natural to generalize the notion of a superposition rule and to consider the mixed superposition rules
that describe the general solution of a system in terms of generic families of particular solutions of (maybe
different) first-order systems and a set of constants. To motivate this concept and to illustrate its
usefulness, we will now provide a series of examples.

Consider a linear differential equation of the form
\begin{equation}\label{simplest}
\frac{dx}{dt}=a(t)x+b(t),
\end{equation}
with $a(t)$ and $b(t)$ being any pair of $t$-dependent functions. It is well known that its general solution
can be written as
\begin{equation}\label{Mix1}
x(t)=x_{(1)}(t)+k x_{(2)}(t),\qquad\qquad k\in\mathbb{R},
\end{equation}
where $x_{(1)}(t)$ and $x_{(2)}(t)$ are generic particular solutions of (\ref{simplest}) and the homogeneous
system
\begin{equation}\label{hom}
\frac{dx}{dt}=a(t)x,
\end{equation}
respectively. In other words, linear systems admit their general solutions to be
described in terms of particular solutions of two different systems, the original one and the homogeneous one,
and a constant.

Let us now turn to Bernoulli equations, i.e. first-order differential equations
\begin{equation}\label{Bernoulli}
\frac{dx}{dt}=a(t)x+b(t)x^n,\qquad\qquad n\neq 1,
\end{equation}
with $a(t)$ and $b(t)$ being two arbitrary $t$-dependent functions. These equations form another class of
first-order differential equations whose general solutions can be obtained from particular solutions of two
different systems and a constant. Indeed, the change of variables $z=x^{1-n}$ transforms (\ref{Bernoulli})
into
$$
\frac{dz}{dt}=(1-n)(a(t)z+b(t)),
$$
whose general solution is $z(t)=z_p(t)+kz_h(t)$, where $k$ is a real constant, $z_p(t)$ is a particular solution of the above
equation, and $z_h(t)$ is a particular solution of
$$
\frac{dz}{dt}=(1-n)a(t)z.
$$
Undoing the previous change of variables, the solution of a Bernoulli equation can be cast in the form
\begin{equation}\label{Mix2}
x(t)=(x_{(1)}^{1-n}(t)+kx_{(2)}^{1-n}(t))^\frac1{1-n},
\end{equation}
where $x_{(1)}(t)$ and $x_{(2)}(t)$ are particular solutions of  (\ref{Bernoulli}) and (\ref{hom}), respectively.

Let us now discuss a family of systems admitting a feature analogue to the above ones. Let $X$ be a Lie
system, so that we can write
\begin{equation}\label{LieMix}
X^t=\sum_{\alpha=1}^rb_\alpha(t)Y_\alpha,
\end{equation}
for certain $t$-dependent functions $b_1(t),\ldots,b_r(t)$ and vector fields $Y_1,\ldots,Y_r$ spanning an
$r$-dimensional real Lie algebra $V$ of vector fields. It is known that there always exists a (local) Lie
group action $\Phi:G\times \mathbb{R}^{n_0}\rightarrow \mathbb{R}^{n_0}$ whose fundamental vector fields
coincide with $V$ \cite[Theorem XI]{Palais}. It can also be proved that the general solution
$x(t)$ of $X$ can be brought into the form
\begin{equation}\label{Mix3}
x(t)=\Phi(g_{(1)}(t);k),
\end{equation}
where $k\in\mathbb{R}^{n_0}$ and $g_{(1)}(t)$ is a particular solution of the Lie system
\begin{equation}\label{LieSysG}
\frac{dg}{dt}=-\sum_{\alpha=1}^rb_\alpha(t)Y^R_{\alpha}(g),\qquad g\in G,
\end{equation}
where each $Y^R_\alpha$ is the unique right-invariant vector field on $G$ corresponding to the fundamental vector field $Y_\alpha\in V$ (see \cite{CGM00,CarRamcinc,Dissertations} for details).
This shows that the general solution of a Lie system can always be described by a particular solution of a
system of the form (\ref{LieSysG}). It is worth noting that this interesting result presents an important
drawback:  frequently the explicit expression of $\Phi$ cannot be determined.

Winternitz-Smorodinsky oscillators \cite{Ev90}-\cite{GPS95} and Milne-Pinney equations \cite{Mil30,Pi50,Ka98,DW10,LA08} can
easily be described by means of systems of first-order differential equations
\begin{equation}\label{WS}
\left\{\begin{aligned}
\frac{dx}{dt}&=p_x,\\
\frac{dp_x}{dt}&=-\omega^2(t)x+\frac{c}{x^3},
\end{aligned}\right.
\end{equation}
with $c$ being a real constant and $x>0$. The general solution $(x(t),p_x(t))$ of any such system can be
determined through two particular solutions $(x_{(1)}(t),p_{(1)}(t))$, $(x_{(2)}(t),p_{(2)}(t))$ of the linear first-order
system
\begin{equation}\label{HarOs}
\left\{\begin{aligned}
\frac{dx}{dt}&=p,\\
\frac{dp}{dt}&=-\omega^2(t)x,
\end{aligned}\right.
\end{equation}
obtained by adding a new variable $p\equiv dx/dt$ to the equation of a $t$-dependent frequency harmonic oscillator,
through the expressions
\begin{equation}\label{Erm}
\begin{aligned}
\!\!\!\!x(t)\!&=\!\frac{\sqrt{2}}{|W|}[k_1x^2_{(1)}(t)+k_2x^2_{(2)}(t)+\sqrt{4k_1k_2-cW^2}x_{(1)}(t)x_{(2)}(t)]^{1/2},\\
\!\!\!\!p(t)\!&=\!\frac{\sqrt{2}[k_1x_{(1)}(t)p_{(1)}(t)\!+\!k_2x_{(2)}(t)p_{(2)}(t)\!+\!\sqrt{k_1k_2\!-\!c(W/2)^2}(p_{(1)}(t)x_{(2)}(t)\!+\!x_{(1)}(t)p_{(2)}(t))]}{|W|[k_1x^2_{(1)}(t)+k_2x^2_{(2)}(t)+\sqrt{4k_1k_2-cW^2}x_{(1)}(t)x_{(2)}(t)]^{1/2}},\\
\end{aligned}
\end{equation}
where $W=x_{(1)}(t)p_{(2)}(t)-p_{(1)}(t)x_{(2)}(t)$ is a constant of motion of system (\ref{HarOs}) and $k_1, k_2$ are two
 real constants (see \cite{CLR08,Mil30,Pi50,Re99}). It is to be remarked that these expressions are interesting because they allow us to study systems appearing in the calculation of invariants for non-quadratic Hamiltonian systems, cosmology, quantum mechanics, Bose-Einstein condensates, etc. \cite{Ka98,DW10,LA08}.

Finally, recall that we aim to introduce a geometric notion covering, as particular cases, expressions of the
type (\ref{Mix1}), (\ref{Mix2}), (\ref{Mix3}), (\ref{Erm}), and standard superposition rules. This leads us to
the following.

\begin{definition} \textbf{(Mixed Superposition Rule)} \label{MSR} A {\it mixed superposition rule} for a system $X$ is a $(m+1)$-tuple $(\Phi,X_{(1)},\ldots,X_{(m)})$ consisting of a function $\Phi:\mathbb{R}^{n_1}\times\ldots\times\mathbb{R}^{n_m}\times \mathbb{R}^{n_0}\rightarrow\mathbb{R}^{n_0}$ and a series of systems $X_{(a)}$ on $\mathbb{R}^{n_a}$, with $a=1,\ldots,m$, such that the general solution $x(t)$ of $X$ can be cast in the form
\begin{equation}\label{MixedSup}
x(t)=\Phi(x_{(1)}(t),\ldots,x_{(m)}(t);k_1,\ldots,k_{n_0}),
\end{equation}
where $x_{(1)}(t),\ldots,x_{(m)}(t)$ is a generic family of particular solutions of $X_{(1)},\ldots,X_{(m)}$,
respectively, and $k_1,\ldots,k_{n_0}$ are constants related to initial conditions.
\end{definition}

\begin{note} Let us stress that the map $\Phi$, describing a mixed superposition rule, does not depend on $t$.
\end{note}

In view of the above definition, a superposition rule for a system $X$ on $\mathbb{R}^{n_0}$ can be viewed
naturally as a mixed superposition rule  of the form
\begin{equation}\label{EspMixSup}
(\Phi:\mathbb{R}^{n_0m}\times\mathbb{R}^{n_0}\rightarrow\mathbb{R}^{n_0},\stackrel{m-{\rm
times}}{\overbrace{X,\ldots,X}}).
\end{equation}

In this language, we can say that linear system (\ref{simplest}) admits the mixed superposition rule given by
$(\Phi_1,X_{(1)},X_{(2)})$, where
$\Phi_1:(x_{(1)},x_{(2)};k)\in\mathbb{R}\times\mathbb{R}\times\mathbb{R}\mapsto x=x_{(1)}+kx_{(2)}\in
\mathbb{R}$, $X_{(1)}=(a(t)x_{(1)}+b(t))\partial/\partial x_{(1)}$, and $X_{(2)}=a(t)x_{(2)}\partial/\partial
x_{(2)}$.

The Bernoulli equations (\ref{Bernoulli}) posses the mixed superposition rule
$(\Phi_2:\mathbb{R}\times\mathbb{R}\times\mathbb{R}\rightarrow\mathbb{R},X_{(1)},X_{(2)})$, where
$\Phi_2(x_{(1)},x_{(2)};k)=(x_{(1)}^{1-n}+kx_{(2)}^{1-n})^\frac{1}{1-n}$,
$X_{(1)}=(a(t)x_{(1)}+b(t)x_{(1)}^n)\partial/\partial x_{(1)}$ and $X_{(2)}=a(t)x_{(2)}\partial/\partial x_{(2)}$.
Meanwhile, every Lie system (\ref{LieMix}) admits a mixed superposition rule $(\Phi;X_{(1)})$, where $X_{(1)}=-\sum_{\alpha=1}^rb_\alpha(t)Y^R_\alpha$ and $\Phi$ is given by (\ref{Mix3}).

Finally, Winternitz-Smorodinsky oscillators (\ref{WS}) or Milne-Pinney equations, written as a first-order
system, admit the mixed superposition rule
$$
\Phi_3:(\xi_{(1)},\xi_{(2)};k_1,k_2)\in {\rm T}^*\mathbb{R}_+\times {\rm T}^*\mathbb{R}_+\times
\mathbb{R}^2\mapsto (x,p)\in {\rm T}^*\mathbb{R}_+,
$$
with $\xi_{(1)}=(x_{(1)},p_{(1)})$, $\xi_{(2)}=(x_{(2)},p_{(2)})$, and
$$
\begin{aligned}
x&=\frac{\sqrt{2}}{|x_{(1)}p_{(2)}-p_{(1)}x_{(2)}|}\left[k_1x^2_{(1)}+k_2x^2_{(2)}+\sqrt{4k_1k_2-c(x_{(1)}p_{(2)}-p_{(1)}x_{(2)})^2}x_{(1)}x_{(2)}\right]^{1/2},\\
p&=\frac{\sqrt{2}\left[k_1x_{(1)}p_{(1)}+k_2x_{(2)}p_{(2)}+\sqrt{k_1k_2-c(p_{(1)}x_{(2)}+p_{(2)}x_{(1)})^2/4}(p_{(1)}x_{(2)}+p_{(2)}x_{(1)})\right]}{|x_{(1)}p_{(2)}-
p_{(1)}x_{(2)}|[k_1x^2_{(1)}+k_2x^2_{(2)}+\sqrt{4k_1k_2-c(x_{(1)}p_{(2)}-p_{(1)}x_{(2)})^2}x_{(1)}x_{(2)}]^{1/2}}.
\end{aligned}
$$
Apart from Milne-Pinney equations, many other systems posses a  similar property, i.e., their general solutions
can be obtained through solutions of a $t$-dependent linear homogeneous system of first-order differential
equation. For instance, Riccati equations, certain second-order Riccati equations,  Kummer-Schwarz equations,
and other systems appearing in the linearization of differential equations admit a similar feature
\cite{PW,Be07,GL99,BR97}.

\section{Mixed superposition rules as flat connections}
Similarly to standard superposition rules, mixed superposition rules can be associated with special flat
connections. This idea, originally developed in \cite{CGM07}, can be applied also in our case.

Consider a mixed superposition rule $(\Phi,X_{(1)},\ldots,X_{(m)})$, with
$\Phi:\mathbb{R}^{n_1}\times\ldots\times\mathbb{R}^{n_m}\times\mathbb{R}^{n_0}\rightarrow\mathbb{R}^{n_0}$,
for a system $X_{(0)}$ on $\mathbb{R}^{n_0}$ (we shall see briefly why this notation is appropriate for our
purposes). The Implicit Function Theorem shows that, fixing a point $p=(x_{(1)},\ldots, x_{(m)})\in
\mathbb{R}^{n_1}\times\ldots\times\mathbb{R}^{n_m}$, the map $\Phi|_p:k\in\mathbb{R}^{n_0}\mapsto
x_{(0)}=\Phi(p;k)\in \mathbb{R}^{n_0}$ can locally be inverted to define a mapping
$\Psi:\mathbb{R}^{n_1}\times\ldots\times\mathbb{R}^{n_m}\times\mathbb{R}^{n_0}\rightarrow\mathbb{R}^{n_0}$
such that
$$
\Psi(x_{(0)},\ldots, x_{(m)})=k,
$$
where $k=(k_1,\ldots,k_{n_0})$ is the only point of $\mathbb{R}^{n_0}$ satisfying
$$
x_{(0)}=\Phi(x_{(1)},\ldots,x_{(m)};k).
$$
Hence, the map $\Psi$ determines an $n_0$-codimensional (generally local) foliation of
$\mathbb{R}^{n_0}\times\ldots\times\mathbb{R}^{n_m}$. As $\Phi$ is a mixed superposition rule, for a generic
family of particular solutions $x_{(0)}(t),\ldots,x_{(m)}(t)$ of $X_{(0)},\ldots,X_{(m)}$, we have
$$
x_{(0)}(t)=\Phi(x_{(1)}(t),\ldots,x_{(m)}(t);k)\Longleftrightarrow \Psi(x_{(0)}(t),\ldots, x_{(m)}(t))=k,
$$
Differentiating the latter expression with respect to $t$, we get
\begin{equation*}
\sum_{a=0}^{m}\sum_{i=1}^{n_a}X_{(a)}^i(t,x_{(a)}(t))\frac{\partial\Psi^j}{\partial x^i_{(a)}}=0,\qquad
j=1,\ldots,n_0,
\end{equation*}
with $\Psi=(\Psi^1,\ldots,\Psi^{n_0})$. Hence,
\begin{equation}\label{Psifi}
\sum_{a=0}^{m}X_{(a)}(t,x_{(a)}(t))\Psi^j(x_{(0)}(t),\ldots,x_{(m)}(t))=Z^t\Psi^j(x_{(0)}(t),\ldots,x_{(m)}(t))=0,
\end{equation}
for $j=1,\ldots,n_0$, where we define $Z$ to be the $t$-dependent vector field on
$\mathbb{R}^{n_0}\times\ldots\times\mathbb{R}^{n_m}$ of the form
\begin{equation}\label{product}
Z=\sum_{a=0}^{m}\sum_{i=1}^{n_a}X_{(a)}^i(t,x_{(a)})\frac{\partial}{\partial x^i_{(a)}}.
\end{equation}
This will be called the {\it direct product}, $Z=X_{(0)}\times\cdots\times X_{(m)}$, of the $t$-dependent vector
fields $X_{(0)},\ldots,X_{(m)}$. More precisely, we have the following definition.
\begin{definition}\label{Ext} Given a set $X_{(0)},\ldots,\!X_{(m)}$ of $t$-dependent vector fields defined, respectively, on $\mathbb{R}^{n_0},\ldots,\mathbb{R}^{n_m}$, their {\it direct product} (or {\it direct prolongation}) $Z=X_{(0)}\times\cdots\times X_{(m)}$ is the unique $t$-dependent vector field $Z$ on $\mathbb{R}^{n_0}\times\ldots\times\mathbb{R}^{n_m}$ such that ${\rm pr}_{a*}Z^t=(X_{(a)})^t$, with ${\rm pr}_a:(x_{(0)},\ldots,x_{(m)})\in\mathbb{R}^{n_0}\times\ldots\times\mathbb{R}^{n_m}\mapsto x_{(a)}\in \mathbb{R}^{n_a},$ for $a=0,\ldots,m$ and $t\in\mathbb{R}$.

We say that a $t$-dependent vector field $Z$ on $\mathbb{R}^{n_0}\times\ldots\times\mathbb{R}^{n_m}$ is a {\it
direct prolongation} if $Z$ can be written as the direct product of a family $t$-dependent vector fields on
the spaces $\mathbb{R}^{n_0},\ldots,\mathbb{R}^{n_m}$. In other words, $Z$ is a direct prolongation if the
vector fields $\{Z^t\}_{t\in\mathbb{R}}$ are projectable onto the spaces $\mathbb{R}^{n_a}$, with
$a=0,\ldots,m$.
\end{definition}

As equalities (\ref{Psifi}) hold for a generic family of particular solutions $x_{(0)}(t),\ldots,x_{(m)}(t)$,
it turns out that a mixed superposition rule $(\Phi,X_{(1)},\ldots,X_{(m)})$ for a system $X_{(0)}$ implies
the existence of $n_0$ common first-integrals, namely $\Psi^1,\ldots,\Psi^{n_0}$, of the vector fields $\{Z^t\}_{t\in\mathbb{R}}$
for the direct product $Z$. These functions give rise to an $n_0$-codimensional foliation $\mathfrak{F}$ such
that the vector fields $Z^t$ are tangent to its leaves, $\mathfrak{F}_k$ with $k\in\mathbb{R}^{n_0}$. Hence,
vector fields from $V^Z$ span an integrable distribution over a dense and open subset of
$\mathbb{R}^{n_0}\times\ldots\times\mathbb{R}^{n_m}$.

The foliation $\mathfrak{F}$ has another important property. Given a level set $\mathfrak{F}_k$ corresponding
to $k=(k_1,\ldots,k_{n_0})$ and $(x_{(1)},\ldots,x_{(m)})\in
\mathbb{R}^{n_1}\times\ldots\times\mathbb{R}^{n_m}$, there is a unique point $x_{(0)}\in \mathbb{R}^{n_0}$
such that  $(x_{(0)},x_{(1)},\ldots,x_{(m)})\in\mathfrak{F}_k$. Then, the projection
\begin{equation}\label{bundle}
{\rm pr}:(x_{(0)},\ldots,x_{(m)})\in \mathbb{R}^{n_0}\times \ldots\times \mathbb{R}^{n_m} \mapsto
(x_{(1)},\ldots,x_{(m)})\in \mathbb{R}^{n_1}\times\ldots\times\mathbb{R}^{n _m}\,,
\end{equation}
induces local diffeomorphisms among the leaves $\mathfrak{F}_k$ of $\mathfrak{F}$ and
$\mathbb{R}^{n_1}\times\ldots\times\mathbb{R}^{n_m}$. Such foliation we will call {\it horizontal} foliation
on the bundle (\ref{bundle}).

This property shows that $\mathfrak{F}$ corresponds to a zero curvature connection $\nabla$ in the bundle
${\rm pr}:\mathbb{R}^{n_0}\times\ldots\times\mathbb{R}^{n_m}\rightarrow
\mathbb{R}^{n_1}\times\ldots\times\mathbb{R}^{n_m}$. Indeed, the restriction of ${\rm pr}$ to each leaf gives
a (local) one-to-one map. In this way, there exists a linear map among vector fields on
$\mathbb{R}^{n_1}\times\ldots\times\mathbb{R}^{n_m}$ and horizontal vector fields tangent to a leaf. This
connection (foliation), along with the vector fields $X_{(1)},\ldots,X_{(m)}$, provides us with a {\it mixed
superposition rule} for $X_{(0)}$ without referring to the map $\Psi$. Indeed, if we take a point $x_{(0)}$
and $m$ particular solutions $x_{(1)}(t),\ldots,x_{(m)}(t)$ of the systems $X_{(1)},\ldots,X_{(m)}$,
respectively, then $x_{(0)}(t)$ is the unique curve in $\mathbb{R}^{n_0}$ such that
$$(x_{(0)}(t),x_{(1)}(t),\ldots, x_{(m)}(t))\in
 \mathbb{R}^{n_0}\times\ldots\times\mathbb{R}^{n_m}
$$
and $(x_{(0)}(0),x_{(1)}(0),\ldots,x_{(m)}(0))$ belong to the same leaf. Thus, it is the foliation
$\mathfrak{F}$ and the systems $X_{(1)},\ldots,X_{(m)},$ what really matters if the {\it mixed superposition
rule} for a system $X_{(0)}$ is concerned.

Conversely, assume that we are given a system $X_{(0)}$, whose general solution we want to analyze. Let
$X_{(1)},\ldots,X_{(m)}$ be the family of systems whose particular solutions will be used to analyze $X_{(0)}$
and $\nabla$ be a flat connection on the bundle ${\rm
pr}:\mathbb{R}^{n_0}\times\ldots\times\mathbb{R}^{n_m}\rightarrow\mathbb{R}^{n_1}\times\ldots\times\mathbb{R}^{n_m}$,
which can  be  integrated  to an $n_0$-codimensional foliation $\mathfrak{F}$ on
$\mathbb{R}^{n_0}\times\ldots\times\mathbb{R}^{n_m}$ such that the vector fields $Z^t$, for $Z=X_{(0)}\times
X_{(1)}\times\cdots\times X_{(m)}$, are tangent to the leaves of $\mathfrak{F}$. Then, the above procedure
provides us with a mixed superposition rule for the system $X_{(0)}$ in terms of solutions of
$X_{(1)},\ldots,X_{(m)}$.

Indeed, let $k\in \mathbb{R}^{n_0}$ enumerate smoothly the leaves $\mathfrak{F}_k$ of $\mathfrak{F}$, i.e.
there exists a smooth map $\iota :\mathbb{R}^{n_0}\rightarrow
\mathbb{R}^{n_0}\times\ldots\times\mathbb{R}^{n_m}$ such that $\iota(\mathbb{R}^{n_0})$ intersects every
$\mathfrak{F}_k$ at a unique point. Then, if $x_{(0)}\in \mathbb{R}^{n_0}$ is the unique point such that
$$(x_{(0)},x_{(1)},\ldots,x_{(m)})\in\mathfrak{F}_k,$$ we can define a {\it mixed superposition rule} for $X_{(0)}$ of the form
$$
x_{(0)}=\Phi(x_{(1)},\ldots,x_{(m)};k),$$ in terms of solutions of the systems $X_{(1)},\ldots,X_{(m)}$.

Let us analyze the above claim in detail. The Implicit Function Theorem shows that there exists a function
$\Psi:\mathbb{R}^{n_0}\times\ldots\times\mathbb{R}^{n_m}\rightarrow\mathbb{R}^{n_0}$ such that
$$\Psi(x_{(0)},\ldots,x_{(m)})=k$$ is equivalent to $(x_{(0)},\ldots,x_{(m)})\in \mathfrak{F}_k$. If we fix
$k$ and take solutions $x_{(1)}(t),\ldots,x_{(m)}(t)$ of $X_{(1)},\ldots,X_{(m)}$, then $x_{(0)}(t)$, defined
by $\Psi(x_{(0)}(t),\ldots,x_{(m)}(t))=k$, is an integral curve of $X_{(0)}$. Indeed, since the vector fields
$Z^t$ are tangent to $\mathfrak{F}$, if $x'_{(0)}(t)$ is a solution of $X_{(0)}$ with the initial value
$x'_{(0)}(0)=x_{(0)}$, then the curve $$t\mapsto (x_{(0)}(t),x_{(1)}(t),\ldots,x_{(m)}(t))$$ lies entirely in
a leaf of $\mathfrak{F}$, so in $\mathfrak{F}_k$. But the point of one leaf is entirely determined by its
projection ${\rm pr}$, so $x'_{(0)}(t)=x_{(0)}(t)$ and $x_{(0)}(t)$ is a solution. Summarizing, we have proved
the following proposition.

\begin{proposition}\label{ConMix}
A mixed superposition rule $(\Phi,X_{(1)},\ldots,X_{(m)})$, with
$\Phi:\mathbb{R}^{n_1}\times\ldots\times\mathbb{R}^{n_m}\times\mathbb{R}^{n_0}\rightarrow\mathbb{R}^{n_0}$,
for a system $X$ on $\mathbb{R}^{n_0}$ amounts to a flat connection (equivalently, horizontal
foliation $\mathfrak{F}$) on the bundle ${\rm
pr}:\mathbb{R}^{n_0}\times\mathbb{R}^{n_1}\times\ldots\times\mathbb{R}^{n_m}
\rightarrow\mathbb{R}^{n_1}\times\ldots\times\mathbb{R}^{n_m}$ such that the vector fields
$\{Z^t\}_{t\in\mathbb{R}}$, associated with the direct product $Z=X\times X_{(1)}\times\ldots\times X_{(m)}$,
are horizontal vector fields with respect to the connection (resp., are tangent to the leaves of
$\mathfrak{F}$).
\end{proposition}
According to the above proposition, the horizontal foliation $\mathfrak{F}$ contains the generalized foliation
$\mathfrak{F}^0$ associated with the generalized distribution $\mathcal{D}$ generated by the Lie algebra
$V^Z={\rm Lie}(\{Z^t\}_{t\in\mathbb{R}})$. One can regard $\mathfrak{F}^0$ as a regular foliation on an open
and dense subset of $\mathbb{R}^{n_0}\times\mathbb{R}^{n_1}\times\ldots\times\mathbb{R}^{n_m}$. If only the
projection
\begin{equation}\label{newprojection} {\rm pr}_{\hat{d}}:\mathbb{R}^{n_0}\times\ldots\times\mathbb{R}^{n_m}\rightarrow\mathbb{R}^{n_0}
\times\ldots\times\widehat{\mathbb{R}^{n_d}}\times\ldots\times\mathbb{R}^{n_m}\,,
\end{equation}
where $\widehat{\mathbb{R}^{n_d}}$ indicates that this space is not included in the direct product, induces
diffeomorphisms on the leaves of $\mathfrak{F}^0$, i.e., it is an injective map on the fibers of
$\mathcal{D}$, then $\mathfrak{F}^0$ can be extended to an $n_d$-codimensional foliation $\mathfrak{F}^{(d)}$, horizontal on the bundle (\ref{newprojection}), which will define a mixed superposition rule for the system $X_{(d)}$. In this way, we get the following.

\begin{proposition}\label{p1} Consider $t$-dependent systems $X_{(a)}$ on $\mathbb{R}^{n_a}$, $a=0,\dots,m$, their direct product $Z=X_{(0)}\times X_{(1)}\times\cdots\times X_{(m)}$, and the corresponding generalized distribution $\mathcal{D}$ generated by the Lie algebra $V^Z={\rm Lie}(\{Z^t\}_{t\in\mathbb{R}})$.
Then, $X_{(d)}$ admits a mixed superposition rule
$$(\Phi',X_{(0)},X_{(1)},\ldots,X_{(d-1)},X_{(d+1)},\ldots,X_{(m)})$$
if and only if the projection (\ref{newprojection}) induces injective maps on the fibers of $\mathcal{D}$.
\end{proposition}

\section{Characterization of systems possessing a mixed superposition rule}

In this section, we first analyze the properties of direct products of $t$-dependent vector fields and other
related notions we define. Subsequently, our results are used to characterize and to analyze systems admitting a mixed
superposition rule.

Let us recall that a $t$-dependent vector field $Z$ on $\mathbb{R}^{n_0}\times\ldots\times\mathbb{R}^{n_m}$ is
called a {\it direct prolongation} if $Z$ can be written as the direct product of a family $t$-dependent
vector fields on the spaces $\mathbb{R}^{n_0},\ldots,\mathbb{R}^{n_m}$, i.e., $Z$ can be then brought into the
form
\begin{equation}\label{product1}
Z=\sum_{a=0}^m\sum_{i=1}^{n_a}Z^i_{a}(t,x_{(a)})\frac{\partial}{\partial x^i_{(a)}},
\end{equation}
for certain functions $Z^{i}_{a}:\mathbb{R}\times\mathbb{R}^{n_a}\rightarrow\mathbb{R}$, with $a=0,\ldots,m$
and $i=1,\ldots,n_a$.

\begin{note}
The term {\it direct prolongation} is coined so as to highlight that this notion is a generalization of the
concept of {\it diagonal prolongations} which appears in the theory of standard superposition rules
\cite{CGM07}.
\end{note}
\noindent In view of (\ref{product1}), the following is obvious.
\begin{lemma}\label{Lie} The Lie bracket of direct prolongations is a direct prolongation.
\end{lemma}

The following lemma provides us with the key property of direct prolongations to characterize systems
admitting mixed superposition rules.
\begin{lemma}\label{Red} Consider a family $Z_1,\ldots,Z_r$ of $t$-{\bf independent} direct prolongations on $\mathbb{R}^{n_0}\times\ldots\times\mathbb{R}^{n_m}$ such that ${\rm pr}_*Z_\alpha$, for $\alpha=1,\ldots,r$, are linearly independent at a generic point. If a vector field $\sum_{\alpha=1}^rf_\alpha Z_\alpha$, with $f_1,\ldots,f_r\in C^\infty(\mathbb{R}^{n_0}\times\ldots\times\mathbb{R}^{n_m})$, is a direct prolongation, then $f_1,\ldots,f_r$ depend only on $x_{(1)},\ldots,x_{(m)}$.
\end{lemma}
\begin{proof} As $Z_1,\ldots,Z_r$ are $t$-independent direct prolongations, we can write
$$
Z_\alpha=\sum_{a=0}^m\sum_{i=1}^{n_a}Z^i_{a\alpha}(x_{(a)})\frac{\partial}{\partial x^i_{(a)}},\qquad
\alpha=1,\ldots,r,
$$
for certain functions $Z^{i}_{a\alpha}:\mathbb{R}^{n_a}\rightarrow\mathbb{R}$. Likewise, if
$\sum_{\alpha=1}^rf_\alpha Z_\alpha$ is a $t$-independent direct prolongation, there
exist functions $B^i_a:\mathbb{R}^{n_a}\rightarrow\mathbb{R}$, with $a=0,\ldots,m$ and $i=1,\ldots,n_a$ such
that
$$
\sum_{\alpha=1}^rf_\alpha Z_\alpha=\sum_{a=0}^m\sum_{i=1}^{n_a}B^i_a(x_{(a)})\frac{\partial}{\partial
x^i_{(a)}}.
$$
Hence,
$$
\sum_{\alpha=1}^rf_\alpha Z^i_{a\alpha}(x_{(a)})=B^i_a(x_{(a)}),\qquad a=0,\ldots,m,\qquad i=1,\ldots,n_a.
$$
In particular, we have the subset of equations
$$
\sum_{\alpha=1}^rf_\alpha Z^i_{a\alpha}(x_{(a)})=B^i_a(x_{(a)}),\qquad a=1,\ldots,m,\qquad i=1,\ldots,n_a.
$$
Since the projections ${\rm pr}_*Z_\alpha$ are linearly independent at a generic point, the above system has a
unique solution $f_1,\ldots,f_r$ whose value is determined by the functions $B^i_a$, $Z^i_{a\alpha}$, for
$a=1,\ldots,m$ and $i=1,\ldots,n_a$, which depend on $x_{(1)},\ldots,x_{(m)}$ only. Hence, $f_1,\ldots,f_r$
depend exclusively on these variables.
\end{proof}

Let us now prove the central result of our paper.

\begin{theorem}\label{MLST}{\bf (The extended Lie-Scheffers Theorem)} A system $X$ admits a mixed superposition rule if and only if it is a Lie system.
\end{theorem}
\begin{proof}
Let
$(\Phi:\mathbb{R}^{n_1}\times\ldots\times\mathbb{R}^{n_m}\times\mathbb{R}^{n_0}\rightarrow\mathbb{R}^{n_0},X_{(1)},\ldots,X_{(m)})$
be a mixed superposition rule for $X$ and let $Z$ be the direct product of $X\times X_{(1)}\times\cdots\times
X_{(m)}$. According to Proposition \ref{ConMix}, the elements of ${\rm Lie}(\{Z^t\}_{t\in\mathbb{R}})$ span a
generalized distribution $\mathcal{D}$ over $\mathbb{R}^{n_0}\times\ldots\times\mathbb{R}^{n_m}$ which is
regular and integrable on an open and dense subset of this space. Further, we can always choose a basis
$Z_1,\ldots,Z_r$ of $\mathcal{D}$ whose elements belong to ${\rm Lie}(\{Z^t\}_{t\in\mathbb{R}})$, so, due to
Lemma \ref{Lie}, are $t$-independent direct prolongations. As $\mathcal{D}$ describes a mixed superposition
rule, the elements of such a basis project,  via ${\rm pr}_*$, onto
$\mathbb{R}^{n_1}\times\ldots\times\mathbb{R}^{n_m}$, giving rise to a family of linearly independent vector
fields at a generic point of this space.

Since each $[Z_\alpha,Z_\beta]$ is a direct prolongation that belongs to $\mathcal{D}$ and $Z_1,\ldots,Z_r$
project to a family of linearly independent vector fields at a generic point of
$\mathbb{R}^{n_1}\times\ldots\times\mathbb{R}^{n_m}$, Lemma \ref{Red} ensures that
$$
[ Z_\alpha, Z_\beta]=\sum_{\gamma=1}^rf_{\alpha\beta\gamma} Z_\gamma,\,\,\,\qquad \alpha,\beta=1,\ldots,r,
$$
for certain $r^3$ functions $f_{\alpha\beta\gamma}$  depending on the variables $x_{(1)},\ldots,x_{(m)}$.
Further, each $[Z_\alpha,Z_\beta]$ is projectable onto $\mathbb{R}^{n_0}$ under the projection ${\rm
pr}_0:(x_{(0)},\ldots,x_{(m)})\in\mathbb{R}^{n_0}\times\ldots\times\mathbb{R}^{n_m}\mapsto x_{(0)}\in
\mathbb{R}^{n_0}$. Hence,
$$
{\rm pr}_{0*}[Z_\alpha,Z_\beta]=[{\rm pr}_{0*}Z_\alpha,{\rm pr}_{0*}
Z_\beta]=[Y_\alpha,Y_\beta]=\sum_{\gamma=1}^rf_{\alpha\beta\gamma}(x_{(1)},\ldots,x_{(m)})Y_\gamma,\quad
\alpha,\beta=1,\ldots,r,
$$
where $Y_\alpha\equiv{\rm pr}_{0*} Z_\alpha$, with $\alpha=1,\ldots,r$, are vector fields on
$\mathbb{R}^{n_0}$. Note that the above equality implies that, for every fixed $(x_{(1)},\ldots,x_{(m)})$, the
vector field $[Y_\alpha,Y_\beta]$ on $\mathbb{R}^{n_0}$ is a linear combination of the vector fields
$Y_1,\ldots,Y_r$ on $\mathbb{R}^{n_0}$, that is, $Y_1,\ldots,Y_r$ span a finite-dimensional real Lie algebra
$V$ of vector fields. Note that this result is due to the fact that all functions  $f_{\alpha\beta\gamma}$
just depend on $x_{(1)},\ldots,x_{(m)}$ only. Now, the direct product $Z=X\times X_{(1)}\times\cdots\times
X_{(m)}$ can be cast in the form
$$
Z=\sum_{\alpha=1}^rb_\alpha Z_\alpha,
$$
where $b_\alpha\in C^{\infty}(\mathbb{R}\times\mathbb{R}^{n_0}\times\ldots\times\mathbb{R}^{n_m})$ for
$\alpha=1,\ldots,r$. Using again the fact that that ${\rm pr}_*Z_1,\ldots,{\rm pr}_* Z_r$ are properly defined
and linearly independent at a generic point, we get out of Lemma \ref{Red} that $b_1,\ldots,b_r$ depend only
on $t$ and $x_{(1)},\ldots,x_{(m)}$. Proceeding as above, we get that ${\rm pr}_{0*} Z^t=X^t$ takes values in
$V$, i.e., $X$ is a Lie system.

Conversely, if $X$ is a Lie system, the Lie-Scheffers Theorem guarantees that it admits a superposition rule.
As superposition rules form a particular class of mixed superposition rules, the theorem follows.
\end{proof}

The extended Lie-Scheffers Theorem not only characterizes systems admitting a mixed superposition rule but
also provides a new tool to ensure that a system is a Lie system: it is enough to find a mixed superposition
rule. For instance, observe that the extended Lie-Scheffers Theorem easily shows that the systems (\ref{WS})
related to Milne-Pinney equations and Winternitz-Smorodinsky oscillators, linear systems of differential equations
(\ref{simplest}), and Bernoulli equations are Lie systems. This easily retrieves as particular cases many of
the results that were obtained in several recent works \cite{CLR08,CL08Diss,CL08,CFLSV} through the standard
Lie-Scheffers Theorem. Additionally, we prove for the first time that all Bernoulli equations are Lie systems, which completes a result merely pointed out in \cite{Dissertations}. In
Section \ref{Appl}, we will describe some new results in this direction.

\begin{corollary} {\bf (extended Lie's condition)} If $X$ is a system admitting a mixed superposition rule $(\Phi:\mathbb{R}^{n_1}\times\ldots\times\mathbb{R}^{n_m}\times\mathbb{R}^{n_0}\rightarrow\mathbb{R}^{n_0},X_{(1)},\ldots,X_{(m)})$, then $X$ admits a Vessiot-Guldberg Lie algebra $V$ such that $\dim\,V\leq \sum_{a=1}^m n_a$.
\end{corollary}
\begin{proof} Following the proof of the extended Lie-Scheffers Theorem, we see that the mixed superposition rule for $X$ induces a family of $t$-independent direct prolongations $Z_1,\ldots, Z_r$ tangent to the leaves of its associated foliation and satisfying that ${\rm pr}_*Z_\alpha$, with $\alpha=1,\ldots,r,$ must be linearly independent at a generic point. Therefore, $r\leq \dim(\mathbb{R}^{n_1}\times\ldots\times\mathbb{R}^{n_m})=\sum_{a=1}^m n_a$. In addition, since the vector fields ${\rm pr}_{0*}Z_\alpha$, with $\alpha=1,\ldots,r,$ span a finite-dimensional real Lie algebra $V$ containing the vector fields $\{X^t\}_{t\in\mathbb{R}}$, our corollary easily follows.
\end{proof}

\begin{note} The above corollary obviously includes the Lie's condition as a particular case. Indeed, it shows that if $X$ admits a superposition rule depending on $m$ particular solutions, then $X$ admits a Vessiot-Guldberg Lie algebra $V$ satisfying that $\dim V\leq n_0\cdot m$.
\end{note}

Let us now discuss our previous results and their relations to the usual notions and properties of
superposition rules.

Observe first that the direct product of $m$ copies of a $t$-dependent vector field $X=\sum_{i=1}^{n_0}X^i(t,x)\partial/\partial x^i$ on $\mathbb{R}^{n_0}$ is a $t$-dependent vector field on
$\mathbb{R}^{n_0m}$ given by (\ref{prolongation}). This is exactly the expression for the diagonal prolongation to $\mathbb{R}^{n_0m}$ of $X$. In other words,
diagonal prolongations are actually direct prolongations involving direct products of several copies of the
same $t$-dependent vector field. On the other hand, some properties of diagonal prolongations are stronger
than those for direct prolongations. In particular, the functions $f_1,\dots,f_r$ that appear in Lemma
\ref{Red} turn out to be just constant for diagonal prolongations (see \cite[Lemma 1]{CGM07}). This simplifies
the proof of the Lie-Scheffers Theorem in comparison with our proof of its extended version.

As mentioned in the introduction, mixed superposition rules provide a new method to study solutions of systems
of differential equations. Although the extended Lie-Scheffers Theorem shows that mixed superposition rules,
like the standard ones, can only be used to study Lie systems, the area of possible applications is much
broader as mixed superposition rules are much more versatile.

First of all, mixed superposition rules can be constructed in various ways and sometimes much easily than the
strict superposition rules. In particular, Lie systems may admit mixed superposition rules in terms of non-Lie
systems. Consider the following trivial example. The equation
$$
\frac{dx}{dt}=0,
$$
is associated with a Lie system (as every autonomous system \cite{CGL08,CGL09,BM09II})  $X=0$. Its general
solution reads $x(t)=k$, with $k\in\mathbb{R}$. Therefore, the map
$\Phi:(x_{(1)};k)\in\mathbb{R}\times\mathbb{R}\mapsto k\in\mathbb{R}$ gives rise to a mixed
superposition rule $(\Phi,X_{(1)})$ for $X$ in terms of any system $X_{(1)}$, e.g. a non-Lie one.

Second, as Proposition \ref{p1} shows, mixed superposition rules may enable us to study general solutions of
different Lie systems at the same time.

\section{Mixed superposition rules in terms of Lie systems}

In spite of the fact that non-Lie systems can appear in mixed superposition rules, most of mixed superposition
rules appearing in the literature involve exclusively Lie systems \cite{Dissertations,CLR08,Be07}. This motivates
the study of this special case. Let us start with some necessary definitions and auxiliary facts to prove the
main results of this section.

\begin{definition} Given a vector space $V$ of vector fields on $\mathbb{R}^{n_0}$, we say that $V$ admits a {\it modular basis} if $V$ possesses a basis of vector fields linearly independent at a generic point of $\mathbb{R}^{n_0}$.
\end{definition}
\begin{note} Note that the term modular basis refers to the fact that the basis of the real Lie algebra $V$ consists of elements which are linearly independent over $C^\infty(\mathbb{R}^{n_0})$ (thus over $\mathbb{R}$). In particular, the space $V$ is finite-dimensional.
\end{note}

\begin{example} The Lie algebra $V=\langle x\partial/\partial x,y\partial/\partial y\rangle$ of vector fields on $\mathbb{R}^2$ obviously possesses a modular basis. On the other hand, the Lie algebra $V=\langle \partial/\partial x,x\partial/\partial x\rangle$ of vector fields on $\mathbb{R}$ does not, as there exist no two vector fields in $V$ linearly independent at a generic (actually any) point of $\mathbb{R}$.
\end{example}

\begin{lemma}\label{Nec} If $V$ is a Lie algebra of vector fields on $\mathbb{R}^{n_0}$
admitting a modular basis, then every basis of $V$ is modular. In particular, if  $Y_1,\ldots,Y_s\in V$ are
vector fields linearly independent over $\mathbb{R}$ and $X\in V$ is of the form
\begin{equation}\label{e1}
X=\sum_{j=1}^sb_j Y_j,
\end{equation}
where $b_1,\ldots,b_s\in C^\infty(\mathbb{R}^{n_0})$, then $b_1,\ldots,b_s$ must be constant.
\end{lemma}
\begin{proof}
Let $\mathcal{D}$ be the generalized distribution spanned by $V$ and let $r=\dim{V}$. Since $V$ admits a
modular basis, the dimension of $\mathcal{D}_p$ is $r$ for a generic point $p$. Suppose that $X_1,\ldots,X_r$
is a basis of $V$ over $\mathbb{R}$. As the vector fields $X_1,\ldots,X_r$ span $\mathcal{D}_p$ for each $p$, the
vectors  $X_1(p),\ldots,X_r(p)$ span an $r$-dimensional space, thus are linearly independent, for a generic
$p$. Hence, $X_1,\ldots,X_r$ is a modular basis of $V$. To prove the last statement, let us observe that we
can enlarge the family $\{Y_1,\ldots,Y_s\}$ to a basis $Y_1,\ldots,Y_r$ of $V$ which is necessarily a modular
basis. We can then write $X$ also in the form
\begin{equation}\label{e2}
X=\sum_{j=1}^rc_j Y_j\,,
\end{equation}
for some constants $c_j$, $j=1,\ldots,r$. As (\ref{e1}) and (\ref{e2}) imply
$$
\sum_{j=1}^s(c_j-b_j)Y_j+\sum_{l=s+1}^rc_lY_l=0
$$
and $Y_1,\ldots,Y_r$ are linearly independent over $C^\infty(\mathbb{R}^{n_0})$, it must be $b_j=c_j$, so
$b_j$ are constants for all $j=1,\ldots,s$.
\end{proof}
\begin{theorem}\label{MT}
If a system $X$ on $\mathbb{R}^{n_0}$ admits a mixed superposition rule in terms of $m$
copies of a Lie system $X_{(1)}$ and $V^{X_{(1)}}_m$ admits a modular basis, then the direct product $Z$ of
$X$ and $m$ times $X_{(1)}$ is a Lie system,  $V^{Z}\simeq V^{X_{(1)}}$ and the projection
$${\rm pr}_0:\mathbb{R}^{n_0}\times \mathbb{R}^{m\cdot n_1}\to\mathbb{R}^{n_0}$$
induces a Lie algebra homomorphism ${\rm pr}_{0*}:V^{Z}\to V^{X}$.
\end{theorem}
\begin{proof}
In view of Proposition \ref{p1}, the vector fields from ${\rm Lie}(\{Z^t\}_{t\in\mathbb{R}})$, where
$$Z=X\times\stackrel{m-{\rm times}}{\overbrace{X_{(1)}\times\ldots\times X_{(1)}}}\,,$$
span a generalized distribution $\mathcal{D}$ over $\mathbb{R}^{n_0}\times \mathbb{R}^{n_1m}$ which is regular
over an open and dense subset of this space and the projection
$${\rm pr}:\mathbb{R}^{n_0}\times \mathbb{R}^{m\cdot n_1}\to\mathbb{R}^{m\cdot n_1}$$
induces injective maps on fibers of $\mathcal{D}$.

We can choose among the elements of $V^Z={\rm Lie}(\{Z^t\}_{t\in\mathbb{R}})$ a basis $Z_1,\ldots,Z_r$ of
$\mathcal{D}$. It follows that such $Z_1,\ldots,Z_r$ are projectable, {\it via} ${\rm pr}_*$, onto
$\mathbb{R}^{n_1m}$ and their projections ${\rm pr }_*Z_1,\ldots,{\rm pr}_*Z_r$ are linearly independent at a
generic point. Note also that such projections are diagonal prolongations, i.e. ${\rm pr}_*Z_\alpha=\widetilde Y_\alpha$ for
certain vector fields $Y_\alpha\in V^{X_{(1)}}$, with $\alpha=1,\ldots,r$. Moreover, we can prove that
$Y_1,\ldots,Y_r$ form a basis of $V^{X_{(1)}}$. Let us start by showing that they span $V^{X_{(1)}}$.
Indeed, as ${\rm pr}_{1*}(V^Z)=V^{X_{(1)}}$, for any element $\bar Y\in V^{X_{(1)}}$ there exists a $t$-independent direct
prolongation $\bar Z\in \mathcal{D}$ such that ${\rm pr}_{1*}\bar Z=\bar Y$. Then, $\bar
Z=\sum_{\alpha=1}^rf_\alpha Z_\alpha$ and ${\rm pr}_*\bar Z=\widetilde{Y}=\sum_{\alpha=1}^rf_\alpha \widetilde
Y_\alpha$, where $f_1,\ldots,f_r$ are certain functions depending only on $x_{(1)},\ldots,x_{(m)}$. Since
$V^{X_{(1)}}_m$ admits a modular basis, Lemma \ref{Nec} ensures that $f_\alpha=c_\alpha$ for certain constants
$c_1,\ldots,c_r$, $\alpha=1,\ldots,r$. Therefore, $\bar Y$ is a linear combination of $Y_1,\ldots,Y_r$.
Additionally, $Y_1,\ldots,Y_r$ must be linearly independent over $\mathbb{R}$ because otherwise $\widetilde
Y_1,\ldots,\widetilde Y_r$ would not be linearly independent at a generic point. Hence, $Y_1,\ldots,Y_r$ form
a basis for $V^{X_{(1)}}$.

Lemma \ref{Red} implies now that
$$
[Z_\alpha,Z_\beta]=\sum_{\gamma=1}^rf_{\alpha\beta\gamma}  Z_{\gamma},\qquad \alpha,\beta=1\ldots,r,
$$
for a unique family of $r^3$ functions $f_{\alpha\beta\gamma}$ depending on the variables
$x_{(1)},\ldots,x_{(m)}$. Hence,
$$
{\rm pr}_*[Z_\alpha,Z_\beta]=[\widetilde Y_\alpha,\widetilde Y_\beta]=\sum_{\gamma=1}^rf_{\alpha\beta\gamma}
\widetilde Y_{\gamma},\qquad \alpha,\beta=1\ldots,r.
$$
As $[\widetilde Y_\alpha,\widetilde Y_\beta]$ belongs to $V^{X_{(1)}}_m$ and $\widetilde Y_1,\ldots,\widetilde
Y_r$ form a modular basis of this Lie algebra, Lemma \ref{Nec} again yields
$f_{\alpha\beta\gamma}=c_{\alpha\beta\gamma}$ for certain $r^3$ constants $c_{\alpha\beta\gamma}$. Hence,
$$
[Z_\alpha,Z_\beta]=\sum_{\gamma=1}^rc_{\alpha\beta\gamma}  Z_{\gamma},\qquad \alpha,\beta=1\ldots,r,
$$
and $Z_1,\ldots,Z_r$ span an $r$-dimensional Lie algebra $V$. Furthermore, as $Y_1,\ldots,Y_r$ share the same
structure constants as $\widetilde Y_1,\ldots,\widetilde Y_r$ and $Z_1,\ldots,Z_r$, then the latter family
of vector fields span a Lie algebra isomorphic to $V^{X_{(1)}}$, i.e. $V^Z\simeq V^{X_{(1)}}$. Note that the
vector fields $\{{\rm pr}_*Z^t\}_{t\in\mathbb{R}}$ are diagonal prolongations of elements $V^{X_{(1)}}$ and
they are therefore spanned by linear combinations (with constant coefficients) of the diagonal prolongations
${\rm pr}_*Z_\alpha=\widetilde Y_\alpha$. In other words, there exists $t$-dependent functions
$b_1(t),\ldots,b_r(t)$ such that ${\rm pr}_*Z^t=\sum_{\alpha=1}^rb_\alpha(t)\widetilde Y_\alpha$, with
$\alpha=1,\ldots,r$. Hence, as $Z^t$ can be described as a unique linear combination of $Z_1,\ldots,Z_r$, it
follows that
$$
Z^t=\sum_{\alpha=1}^rb_\alpha(t)Z_\alpha.
$$
In other words, $Z$ is a Lie system related to a Vessiot-Guldberg Lie algebra $V^Z$. Obviously $V^Z\simeq
V^{X_{(1)}}$ and, as ${\rm pr}_{0*}{\rm Lie}(\{Z^t\}_{t\in\mathbb{R}})={\rm Lie}(\{{\rm
pr}_{0*}Z^t\}_{t\in\mathbb{R}})={\rm Lie}(\{X^t\}_{t\in\mathbb{R}})=V^{X}$, we have $V^{X}= {\rm
pr}_{0*}V^{Z}$.
\end{proof}

Let us given a second result which allows us to link, under certain conditions, the properties of a Lie system to those of the Lie systems involved in its mixed superposition rules.

\begin{theorem}\label{MT2}
If a system $X$ on $\mathbb{R}^{n_0}$ admits a mixed superposition rule
$$
(\Phi:\mathbb{R}^{n_1}\times\ldots\times\mathbb{R}^{n_m}\times\mathbb{R}^{n_0}\rightarrow\mathbb{R}^{n_0},X_{(1)},\ldots,X_{(m)})
$$in terms of $m$ Lie systems $X_{(1)},\ldots,X_{(m)}$, and, for a certain $k\in\mathbb{R}^{n_0}$, the mapping
$$\Phi_k:(x_{(1)},\ldots,x_{(m)})\in\mathbb{R}^{n_1}\times\ldots\times\mathbb{R}^{n_m}\mapsto \Phi(x_{(1)},\ldots,x_{(m)};k)\in\mathbb{R}^{n_0},\qquad k\in \mathbb{R}^{n_0},$$
has open and dense image in $\mathbb{R}^{n_0}$, then its tangent map induces a Lie algebra epimorphism $\Phi_{k*}:V^{\widehat X}\to V^{X}$, with $\widehat X=X_{(1)}\times\ldots\times X_{(m)}$.
\end{theorem}
\begin{proof}
Note that, under the above assumptions,
\begin{equation*}
x(t)=\Phi(x_{(1)}(t),\ldots,x_{(m)}(t);k)=\Phi_k(x_{(1)}(t),\ldots,x_{(m)}(t)),
\end{equation*}
is a solution of $X$ for every generic family of particular solutions $x_{(1)}(t),\ldots,x_{(m)}(t)$ of the systems $X_{(1)},\ldots,X_{(m)}$, respectively. Differentiating in terms of $t$ and using that ${\rm Im}\, \Phi_{k}$ is dense and open in $\mathbb{R}^{n_0}$, we obtain
$$
X^t=\Phi_{k*}(X_{(1)}\times\ldots\times X_{(m)})_t=\Phi_{k*}\widehat X^t,\qquad\forall t\in\mathbb{R}.
$$
As $X_{(1)},\ldots,X_{(m)}$ are Lie systems, the vector fields corresponding to
$$\begin{array}{ccccccccc}
X_{(1)}&\times &0&\times&\ldots&\times& 0&\times& 0,\\
0&\times &X_{(2)}&\times& \ldots&\times& 0&\times& 0, \\
& &&& \ldots&& && \\
0&\times &\ldots&\times& \ldots&\times& 0&\times& X_{(m)}
\end{array}
$$
span a finite-dimensional Lie algebra of vector fields containing $V^{\widehat X}$, which shows that $\widehat X$ is a Lie system. From here, it can easily be proved that $\Phi_{k*}$ is well defined over $V^{\widehat X}$, i.e. its elements are projectable under $\Phi_{k*}$, and  $V^X=\Phi_{k*}V^{\widehat X}$.
\end{proof}

\begin{corollary}\label{PracRule} If a non-zero Lie system $X$ admits a mixed superposition rule $(\Phi:\mathbb{R}^{n_1m}\times\mathbb{R}^{n_0}\rightarrow\mathbb{R}^{n_0},X_{(1)},\ldots,X_{(1)})$ in terms of $m$ particular solutions of a Lie system $X_{(1)}$ such that $V^{X_{(1)}}$ is a simple Lie algebra and at least one of the following conditions holds:
\begin{enumerate}
 \item $V^{X_{(1)}}_m$ admits a modular basis,
\item There exists a $k\in\mathbb{R}^{n_0}$ such that ${\rm Im}\,\Phi_k$ is open and dense in $\mathbb{R}^{n_0}$,
\end{enumerate}
then $V^X\simeq V^{X_{(1)}}$.
\end{corollary}
\begin{proof} If condition $(1)$ holds, Theorem \ref{MT} shows that $V^Z\simeq V^{X_{(1)}}$, with $Z\!=\!X\times\! \stackrel{m-{\rm times}}{\overbrace{X_{(1)}\!\times\!\ldots\!\times\! X_{(1)}}}$. So if $V^{X_{(1)}}$ is simple, then $V^Z$ is simple and any non-trivial Lie algebra homomorphism from $V^{Z}$ is an isomorphism. In particular, $V^X= {\rm pr}_*(V^Z)\simeq V^Z\simeq V^{X_{(1)}}$.

In a similar way, if $(2)$ is satisfied, Theorem \ref{MT2} states that there exists a Lie algebra epimorphism from $V^{\widehat X}$ to $V^X$, with $\widehat X= X_{(1)}\times\ldots\times X_{(1)}$ ($m$-times). As $V^{X_{(1)}}$ is simple, $V^{\widehat X}\simeq V^{X_{(1)}}$ is also. Hence, as $V^X\neq \{0\}$, then $\Phi_{k*}:V^{\widehat X}\rightarrow V^X$ is an isomorphism and $V^X\simeq V^{\widehat X}\simeq V^{X_{(1)}}$.
\end{proof}

Observe that Theorems \ref{MT} and \ref{MT2} as well as Corollary \ref{PracRule} provide information about the Vessiot-Guldberg Lie
algebras associated to a system admitting a mixed superposition rule in terms of other Lie systems. This will be used in the following
section to describe a new infinite family of Lie systems with interesting applications as well as to give
simple proofs of some known results about Kummer-Schwarz, Ermakov, Riccati and second-order Riccati equations.

\section{Mixed superposition rules and the Riccati hierarchy}\label{Appl}

Consider a linear homogeneous differential equation
\begin{equation}\label{linear}
\frac{d^sx}{dt^s}=-\sum_{l=0}^{s-1}b_l(t)\frac{d^{l}x}{dt^l},
\end{equation}
with $b_0(t),\ldots,b_{s-1}(t)$ being arbitrary functions of time, $d^0x/dt^0\equiv x$ and $s\geq 2$. This
equation is invariant under dilations. This symmetry induces a change of variables $y=x^{-1}dx/dt$, which
transforms the above linear system into a differential equation
\begin{equation}\label{HRiccatiE}
\frac{d^{s-1}y}{dt^{s-1}}=F_{b}\left(t,y,\frac{dy}{dt},\ldots,\frac{d^{s-2}y}{dt^{s-2}}\right),
\end{equation}
for a nonlinear function $F_{\it b}:\mathbb{R}\times\mathbb{R}^{s-1}\rightarrow\mathbb{R}$ whose form depends on $b=(b_0(t),\ldots,b_{s-1}(t))$.

Equations (\ref{HRiccatiE}) are called {\it higher-order Riccati equations}. These equations describe almost the whole family of differential equations of the so-called {\it Riccati hierarchy} \cite{GL99,CRS05}, which is of importance in the study of soliton solutions for PDEs and other relevant physical topics \cite{CL10SecOrd2,GL99}.  Its first members read \cite{GL99}
\begin{eqnarray}
&\dfrac{dy}{dt}=-b_0(t)-b_1(t)y-y^2,\qquad s=2,\\
&\dfrac{d^2y}{dt^2}=-3y\dfrac{dy}{dt}-y^3-b_0(t)-b_1(t)y-b_2(t)\left(y^2+\dfrac{dy}{dt}\right),\qquad s=3. \label{RicSOr}
\end{eqnarray}
The first is a well-known {\it Riccati equation}, which is almost ubiquitous in physics and mathematics
\cite{Dissertations}. The second one can be understood as a generalization of the Painlev\'e-Ince equation
\cite{In65,KL09} and it has been recently frequently studied \cite{CL10SecOrd2,CRS05}.

The general solution
$y(t)$ of (\ref{HRiccatiE}) can be written as
\begin{equation}\label{iniMix}
y(t)=\varphi\left(\!x_{(1)}(t),\ldots,x_{(s)}(t),\!\frac{dx_{(1)}}{dt}(t),\ldots,\frac{dx_{(s)}}{dt}(t),k\!\right)\!\!\equiv\!\! \left(\sum_{a=1}^s\!k_a\frac{dx_{(a)}}{dt}(t)\right)\!\!\left(\sum_{a=1}^s\!k_ax_{(a)}(t)\!\!\right)^{\!\!-1}
\end{equation}
in terms of a generic family $x_{(1)}(t),\ldots,x_{(s)}(t)$ of solutions of (\ref{linear}) and a generic
 $k=(k_1,\ldots,k_s)\in\mathbb{R}^s$. Note that the particular solution induced by $k$ and $\lambda k$ is the same for every $\lambda\in\mathbb{R}/\{0\}$, i.e.
$$\varphi\!\left(\!x_{(1)}(t),\ldots,x_{(s)}(t),\!\frac{dx_{(1)}}{dt}(t),\ldots,\!\frac{dx_{(s)}}{dt}(t),k\!\right)\!=\!\varphi\!\left(\!x_{(1)}(t),\ldots,x_{(s)}(t),\! \frac{dx_{(1)}}{dt}(t),\ldots,\!\frac{dx_{(s)}}{dt}(t),\lambda k\!\right).$$
Hence, expression (\ref{iniMix}) enables us to describe a generic solution $y(t)$ of (\ref{HRiccatiE}) in terms of $x_{(1)}(t),\ldots,x_{(m)}(t)$, arbitrary constants $k_1,\ldots,k_{s-1}$, and $k_s=1$. This implies that the successive derivatives of $y(t)$ can be brought into the
form
\begin{equation}\label{MixSupHOR}
\frac{d^ly}{dt^l}(t)=\Phi^l\left(x_{(1)}(t),\ldots,x_{(s)}(t),\ldots,\frac{d^{s-1}x_{(1)}}{dt^{s-1}}(t),\ldots,\frac{d^{s-1}x_{(s)}}{dt^{s-1}}(t);k_1,\ldots,k_{s-1}\right),
\end{equation}
for certain functions $\Phi^l:\mathbb{R}^{s^2}\times\mathbb{R}^{s-1}\rightarrow\mathbb{R}$, with $l=1,\ldots,s-2$.
Let us view these relations as a mixed superposition rule.

The higher-order differential equations (\ref{linear}) and (\ref{HRiccatiE}) can be written as systems of
first-order differential equations
\begin{equation}\label{Syst2}\left\{
\begin{aligned}
\frac{du^i}{dt}&=u^{i+1},\qquad\qquad\qquad\qquad i=0,\ldots,s-2,\\
\frac{du^{s-1}}{dt}&=-\sum_{l=0}^{s-1}b_l(t)u^{l},
\end{aligned}\right.
\end{equation}
and
\begin{equation}\label{Syst3}
\left\{
\begin{aligned}
\frac{dv^i}{dt}&=v^{i+1},\qquad\qquad\qquad\qquad i=0,\ldots,s-3,\\
\frac{dv^{s-2}}{dt}&=F_{b}(t,v^0,\ldots,v^{s-2}),
\end{aligned}\right.
\end{equation}
where $u^0=x$ and $v^0=y$. Let $X^{b}_{(1)}$ and $X_{F_{b}}$ be the $t$-dependent vector fields associated with (\ref{Syst2}) and (\ref{Syst3}), respectively.  Expressions (\ref{iniMix}) and (\ref{MixSupHOR}) allow
us to write the general solution ${\bf v}(t)=(v^0(t),v^1(t),\ldots,v^{s-2}(t))$ of $X_{F_{b}}$ in terms of a generic
family of particular solutions ${\bf u_{(a)}}(t)=(u^0_{(a)}(t),\ldots,u^{s-1}_{(a)}(t))$, with $a=1,\ldots,s$,
of the system $X^{b}_{(1)}$ and constants $k_1,\ldots,k_{s-1}$. More geometrically, if we define
$$
\Phi^0({\bf u}_{(1)},\ldots,{\bf
u}_{(s)};k_1,\ldots,k_{s-1})=\left(\sum_{a=1}^{s-1}k_au^1_{(a)}+u^1_{(s)}\right)\left(\sum_{a=1}^{s-1}k_au^0_{(a)}+u^0_{(s)}\right)^{-1},
$$we can construct a map
$$
\Phi:({\bf u}_{(1)},\ldots,{\bf u}_{(s)};k_1,\ldots,k_{s-1})\in \mathbb{R}^{s^2}\times\mathbb{R}^{s-1}\mapsto
\Phi({\bf u}_{(1)},\ldots,{\bf u}_{(s)};k_1,\ldots,k_{s-1})\in \mathbb{R}^{s-1},
$$
with $\Phi=(\Phi^0,\ldots,\Phi^{s-2}),$ which enables us to write the general solution ${\bf v}(t)$ of $X^{F_{b}}$ in
terms of $s$ generic particular solutions of $X^{b}_{(1)}$ as
$$
{\bf v}(t)=\Phi({\bf u}_{(1)}(t),\ldots,{\bf u}_{(s)}(t);k_1,\ldots,k_{s-1}).
$$
That is, we have defined a mixed superposition rule $(\Phi,\stackrel{s-{\rm
times}}{\overbrace{X^{b}_{(1)},\ldots,X^{b}_{(1)}}})$ for $X_{F_{b}}$. Therefore, in view of the extended Lie Theorem, each $X_{F_b}$ is a Lie system. For $s=1$ and $s=2$, this retrieves that Riccati equations and second-order Riccati equations of the form (\ref{RicSOr}) are Lie systems \cite{PW,CL10SecOrd2} and shows that these cases are particular instances of a general property of the systems $X_{F_{b}}$ related to the  members of Riccati hierarchy (\ref{HRiccatiE}).

 Depending on the particular form of the functions $b_0(t),\ldots,b_{s-1}(t)$, the minimal Lie algebras for systems (\ref{Syst3}) range from a unidimensional one, when $b_0(t),\cdots,b_{s-1}(t)$ are constant, to the Lie algebra $V^{X_{F_{b}}}$ corresponding to the case $b=(b_0(t),\ldots,b_{s-1}(t))$ where the vectors $(b_0(t),\ldots,b_{s-1}(t))\in\mathbb{R}^s$, with $t\in\mathbb{R}$, span the whole $\mathbb{R}^s$. We will hereafter focus on this latter case, as it is a generic case and whose $V^{X_{F_{b}}}$ describes a Vessiot-Guldberg Lie algebra for all the systems (\ref{Syst3}) with the same $s$. Note that, from the change of variables mapping (\ref{linear}) to (\ref{HRiccatiE}), it easily follows that $\dim\, V^{X_{F_{b}}}>1$ (cf. \cite{GL99,CRS05}).

As every $X^{b}_{(1)}$ is a linear system, it is a Lie system.
Moreover, it is easy to prove that
$V^{X^{b}_{(1)}}$ is spanned by the linear vector fields
\begin{equation}\label{Fam}
X_{i,j}=u^j\frac{\partial}{\partial u^i},\quad i,j=0,\ldots,s-1.\end{equation}
Indeed, the linear spaces spanned by $\{(X^{b}_{(1)})^t\}_{t\in\mathbb{R}}$, on one hand, and the vector fields
$X_{s-1,j}$, with $j=0,\ldots,s-1$, and $\Delta=X_{0,1}+\ldots+X_{s-2,s-1}+X_{s-1,0}$, on the other, are the same. So, they are also the smallest Lie algebras containing their elements. Since
$$[X_{i,j},\Delta]=X_{i-1,j}-X_{i,j+1},\quad [X_{i,s-1},\Delta]=X_{i-1,s-1}-X_{i,0}, \quad i=1,\ldots,s-1,\quad j=0,\ldots,s-2,
$$
it inductively follows that the successive Lie brackets of  elements of
$\{(X^{b}_{(1)})^t\}_{t\in\mathbb{R}}$ span the vector fields $X_{i,j}$, which clearly generate $V^{X^{b}_{(1)}}$. The elements of this family form a basis of a real Lie algebra isomorphic to $\mathfrak{gl}(s,\mathbb{R})$. Moreover, it is straightforward to prove that their prolongations to $\mathbb{R}^{s^2}$ form a modular basis of a real Lie algebra isomorphic to $\mathfrak{gl}(s,\mathbb{R})$. Hence, Theorem \ref{MT} ensures that $V^{X^{b}_{(1)}}_s\simeq V^{X^{b}_{(1)}}\simeq \mathfrak{gl}(s,\mathbb{R})$ is a Lie extension of $V^{X_{F_{b}}}$.  As $\mathfrak{gl}(s,\mathbb{R})$ admits a {\it Levi decomposition} $\mathfrak{sl}(s,\mathbb{R})\oplus \mathbb{R}$, the Lie algebra $V^{X^{b}_{(1)}}$ admits only proper ideals isomorphic to $\mathbb{R}$ and $\mathfrak{sl}(s,\mathbb{R})$, respectively. Then, as $\dim\, V^{X_{F_{b}}}>1$, we obtain that $V^{X_{F_{b}}}$ must be isomorphic to $\mathfrak{sl}(s,\mathbb{R})$ or  $\mathfrak{gl}(s,\mathbb{R})$.

There exists no real Lie algebras of vector fields isomorphic to $\mathfrak{gl}(2,\mathbb{R})$ and $\mathfrak{gl}(3,\mathbb{R})$ over $\mathbb{R}$ and $\mathbb{R}^2$, respectively (cf. \cite{Lie1880,GLKP90}). Hence, the above result, for $s=1$ and $s=2$, shows that $V^{X_{F_{b}}}$ is isomorphic to $\mathfrak{sl}(2,\mathbb{R})$ and $\mathfrak{sl}(3,\mathbb{R})$, correspondingly. Therefore, the Riccati equations and second-order Riccati
equations (\ref{RicSOr}) can be described as Lie systems related to Vessiot-Guldberg Lie algebras isomorphic to
$\mathfrak{sl}(2,\mathbb{R})$ and $\mathfrak{sl}(3,\mathbb{R})$. This has been proved by a direct application
of the Lie-Scheffers Theorem in \cite{Dissertations,CL10SecOrd2} that had required the determination of the
corresponding complicated Lie algebras of non-linear vector fields. Here, this result appears as a geometric consequence of admitting a mixed superposition rule that requires no calculation. Actually, we can determine a Vessiot-Guldberg algebra for every $X_{F_b}$ with a generic $s$ as follows.

\begin{proposition} The system $X_{F_{b}}$ associated with a Riccati equation of order $s$ is a Lie system possessing a Lie-Vessiot algebra isomorphic to $\mathfrak{sl}(s,\mathbb{R})$.
\end{proposition}
\begin{proof} As $X_{F_{b}}$ is a Lie system, $V^{X_{F_{b}}}$ is a Vessiot-Guldberg Lie algebra for $X_{F_{b}}$. Let us analyze its algebraic structure through Theorem \ref{MT2}.

It is can readily be proved that the mixed superposition rule for $X^{F_{b}}$ can be brought into the form
\begin{equation*}
\Phi(p;k)=\left(\Phi^0(p;k),D_T\Phi^0(p;k),\ldots,D^{s-2}_T\Phi^0(p;k)\right)\in\mathbb{R}^{s-1},
\end{equation*}
where $k\in\mathbb{R}^{s-1}$ and we denote $p=({\bf u_{(1)}},\ldots,{\bf u_{(s)}})\in\mathbb{R}^{s^2}$ and 
$$
D_T=\sum_{a=1}^s\sum_{i=0}^{s-2}u^{i+1}_{(a)}\frac{\partial}{\partial u^i_{(a)}}.
$$
From the above expressions, it can easily be demonstrated that $\Phi_k(\cdot)=\Phi(\cdot;k)$ is surjective for a generic $k$. Therefore, Theorem \ref{MT2} implies that $\Phi$ induces a Lie algebra epimorphism $\Phi_{k*}:V_s^{X^{b}_{(1)}}\rightarrow V^{X_{F_{b}}}$ for a generic $k\in\mathbb{R}^{s-1}$. Let us prove that this map has a nontrivial kernel.

The function $\Phi^0_k(\cdot)\equiv \Phi^0(\cdot\,;k)$ is homogeneous, i.e. $\Phi^0_k(\lambda{\bf u_{(1)}},\ldots,\lambda{\bf u_{(s)}})=\Phi^0_k({\bf u_{(1)}},\ldots,{\bf u_{(s)}})$ for all $\lambda\in\mathbb{R}/\{0\}.$ From here, we obtain that the functions $D^j_T\Phi^0_k$, with $j=1,\ldots,s-1$, are also homogeneous. Consequently, $\Phi_k$ is homogeneous.  Since the flow of the vector field $$X_0=\sum_{i=0}^{s-1}\sum_{a=1}^{s}u^i_{(a)}\partial/\partial u^i_{(a)}\in V_s^{X^{b}_{(1)}}$$ reads $$g:(\lambda;{\bf u_{(1)}},\ldots,{\bf u_{(s)}})\in\mathbb{R}\times\mathbb{R}^{s^2}\mapsto e^\lambda({\bf u_{(1)}},\ldots,{\bf u_{(s)}})\in\mathbb{R}^{s^2}\,,$$
then
$$
(\Phi_{k*}X_0)({\bf u_{(1)}},\ldots,{\bf u_{(s)}})=\frac{d}{d\lambda}\bigg|_{\lambda=0}\Phi_k\circ g(\lambda;{\bf u_{(1)}},\ldots,{\bf u_{(s)}})=\frac{d}{d\lambda}\bigg|_{\lambda=0}\Phi_k({\bf u_{(1)}},\ldots,{\bf u_{(s)}})=0
$$
and $\Phi_{k*}$ has a nontrivial kernel. Taking into account that $V_s^{X^{b}_{(1)}}\simeq \mathfrak{sl}(s,\mathbb{R})\oplus\mathbb{R}$ and $\dim(V^{X_{F_{b}}})> 1$, we see that $V^{X_{F_{b}}}\simeq \mathfrak{sl}(s,\mathbb{R})$.
\end{proof}

Note that most of our above procedure does not depend on the explicit form of the mixed superposition rule. Hence, it can be applied, under slight modifications, to analyze systems admitting a mixed superposition rule in terms of solutions of a linear system $X^b_{(1)}$. This occurs, for instance, in the study of second-order Kummer-Schwarz equations and Milne-Pinney equations \cite{Dissertations,Be07,BR97}. In these cases, we have
$$
X_F=v\frac{\partial}{\partial x}+F(t,x,v)\frac{\partial}{\partial v}\,,\qquad X^b_{(1)}=v\frac{\partial}{\partial x}-\omega^2(t)x\frac{\partial}{\partial v}\,,
$$
where $F$ is a $t$-dependent function, related to one of the previous second-order differential equations, and $X^b_{(1)}$ is associated with a $t$-dependent frequency harmonic oscillator $d^2x/dt^2=-\omega^2(t)x$. Thus, $b_0(t)=\omega^2(t)$, $b_{1}(t)=0$, and $V_2^{X^b_{(1)}}\simeq \mathfrak{sl}(2,\mathbb{R})$ admits a modular basis as before. As $V^{X_F}\neq 0$ and $V^{X^b_{(1)}}$ is simple, Corollary \ref{PracRule} shows that the system $X_F$, corresponding to a Kummer-Schwarz or Milne-Pinney equation, is a Lie
system possessing a Vessiot-Guldberg Lie algebra isomorphic to $\mathfrak{sl}(2,\mathbb{R})$. This was recently discovered  through the Lie-Scheffers Theorem \cite{PW,CLR08,CGL11}, but our method yields this result much simplier and in a unified form.

\section{Conclusions and Outlook}
Following an idea briefly suggested in \cite{CGM07}, we have proposed a definition of a mixed superposition rule
that generalizes the concept of a superposition rule of Lie and Scheffers. We have characterized systems
admitting a mixed superposition rule as certain flat connections and proved a result, called here the extended
Lie-Scheffers Theorem, stating that only Lie systems admit mixed superposition rules. Additionally, we have
shown that mixed superposition rules are more versatile than the standard ones and may be utilized to analyze
simultaneously the general solutions of different Lie systems. The extended Lie-Scheffers Theorem provided us
also with a new powerful tool for recognizing Lie systems which was then applied for retrieving, in a simple
way, some results known from the recent literature.

Our methods have been illustrated by various examples of physical and mathematical interest, in particular, by
relevant results about the hierarchy of Riccati equations. This gave rise to the description of a new
interesting infinite family of Lie systems with relevant applications. Furthermore, our methods seem to apply to certain differential equations  appearing in the linearization of higher-order differential equations \cite{Be07,Da97}.

Finally, it is natural to search further for a characterization of families of systems of first-order differential equations admitting a $t$-dependent common mixed superposition rule  \cite{CGL09}. This could lead to the joint analysis of the whole Riccati hierarchy as well as other systems of differential equations possessing common $t$-dependent mixed superposition rules. We aim to study these questions and other possible applications in a future work.

\section*{Acknowledgements}

J. de Lucas acknowledges partial financial support by the research projects MTM2010-12116-E,
FMI43/10 (DGA) and E24/1 (DGA). Research of J. Grabowski is financed by the Polish Ministry of Science and
Higher Education under the grant N N201 416839.


\begin{thebibliography}{19}
\bibitem{LS}   S. Lie and G. Scheffers,
{\sl Vorlesungen \"uber continuierliche Gruppen mit geometrischen
 und anderen Anwendungen},
  Teubner, Leipzig, 1893.

\bibitem{PW}
P. Winternitz, {\it Lie groups and solutions of nonlinear differential equations}, in: {\sl Nonlinear
Phenomena}, Lecture Notes in Phys. {\bf 189}, Springer-Verlag, Berlin, 1983,  263--305.

\bibitem{CGM00}
J.F. Cari\~nena, J. Grabowski and G. Marmo, {\sl Lie-Scheffers systems: a geometric approach}, Bibliopolis, Naples, 2000.

\bibitem{CGM07}
J.F. Cari\~nena,  J. Grabowski and G. Marmo, {\it Superposition rules, lie theorem and partial differential
equations}, Rept. Math. Phys. {\bf 60}, 237--258 (2007).

\bibitem{CarRamcinc}
J.F. Cari\~nena and  A. Ramos {\it A new geometric approach to Lie systems and physical applications}, Acta
Appl. Math. {\bf 70}, 43--69 (2002).

\bibitem{Dissertations}
J. F. Cari\~nena and J. de Lucas, {\sl Lie systems: theory, generalizations, and applications},
Dissertationes Math. {\bf 479}, Institute of Mathematics of the Polish Academy of Sciences, Warsaw, 2011.

\bibitem{FLV10}
R. Flores-Espinoza, J. de Lucas and Y.M. Vorobiev, {\it Phase splitting for periodic Lie systems}, J. Phys. A
{\bf 43}, 205208  (2010).

\bibitem{Fl10}
R. Flores-Espinoza,
{\it Periodic first integrals for Hamiltonian systems of Lie type},
Int. J. Geom. Methods Mod. Phys. {\bf 8}, 1169--1177 (2011).

\bibitem{Clem06}
J. Clemente-Gallardo,
{\it On the relations between control systems and Lie systems}, in: {\sl Groups, geometry and physics}, Monogr. Real Acad. Ci. Exact. F\'is.-Qu\'im. Nat. Zaragoza vol. 29, Acad. Cienc. Exact. F\'is. Qu\'im. Nat. Zaragoza, Zaragoza, 2006, 65--78.

\bibitem{CLR08}
J.F. Cari\~nena, J. de Lucas, and M.F. Ra\~nada, {\it Recent applications of the theory of Lie systems in
Ermakov systems}, SIGMA Symmetry Integrability Geom. Methods Appl. {\bf 4}, 031 (2008).


\bibitem{CL08Diss}  J. F. Cari\~nena and  J. de Lucas,
{\it Applications of Lie systems in dissipative Milne--Pinney equations}, Int. J. Geom. Methods Mod. Phys.
{\bf 6}, 683--699 (2009).

\bibitem{HWA83}
J. Harnad, P. Winternitz and R.L. Anderson, {\it Superposition principles for matrix Riccati equations}, J.
Math. Phys. {\bf 24}, 1062--1072 (1983).

\bibitem{CL08}
J. F. Cari\~nena and J. de Lucas, {\it A nonlinear superposition rule for solutions of the Milne--Pinney
equation}, Phys. Lett. A {\bf 372}, 5385--5389 (2008).


\bibitem{CL10SecOrd2}
J.F. Cari\~nena and J. de Lucas, {\it Superposition rules and second-order Riccati equations}, J. Geom. Mech.
{\bf 3}, 1--22 (2011).

\bibitem{CGL11}
J.F. Cari\~nena, J. Grabowski and J. de Lucas, {\it Superposition rules for higher-order systems
and their applications}, arXiv:1111.4070.

\bibitem{JP09}
J.A. L\'azaro-Cam\'i and J.P. Ortega, {\it Superposition rules and stochastic Lie-Scheffers systems}, Ann.
Inst. Henri Poincar\'e Probab. Stat. {\bf 45}, 910--931 (2009).

\bibitem{CGL08}
J.F. Cari\~nena, J. Grabowski and J. de Lucas, {\it Quasi-Lie schemes: theory and applications}, J. Phys. A
{\bf 42}, 335206 (2009).

\bibitem{CGL09}
J.F. Cari\~nena, J. Grabowski and J. de Lucas, {\it Lie families: theory and applications}, J. Phys. A {\bf
43}, 305201 (2010).

\bibitem{Ve95}
M.E. Vessiot, {\it Sur quelques \'equations diff\'erentielles ordinaires du second ordre}, Ann. Fac. Sci. Toulousse {\bf 3}, F1--F26  (1895).

\bibitem{BecGagHusWin90}
J. Beckers, L. Gagnon, V. Hussin and P. Winternitz, {\it Superposition formulas for nonlinear superequations},
J. Math. Phys. {\bf 31}, 2528--2534 (1990).

\bibitem{In86}
E.L. Ince, {\sl Ordinary Differential Equations}, Dover Publications, New York, 1944.

\bibitem{WintSecond}
C. Rogers, W.K. Schief and P. Winternitz, {\it Lie-theoretical generalization and discretization of the
Pinney Equation}, J. Math. Anal. Appl. {\bf 216}, 246--264 (1997).


\bibitem{In65}
A. Inselberg, {\it On classification and superposition principles for nonlinear operators}, thesis (Ph.D.) -
University of Illinois at Urbana-Champaign. ProQuest LLC, Ann Arbor, MI, 1965.

\bibitem{Be07}
L.M. Berkovich, {\it Method of factorization of ordinary differential operators and some of its applications},
Appl. Anal. Discrete Math. {\bf 1}, 122--149 (2007).

\bibitem{Mil30}
W.E. Milne, {\it The numerical determination of characteristic numbers}, Phys. Rev. {\bf 35}, 863--867 (1930).


\bibitem{GGG11}
P. Guha, A. Ghose Choudhury and B. Grammaticos, {\it Dynamical studies of equations from the Gambier Family},
SIGMA Symmetry Integrability Geom. Methods Appl. {\bf 7}, 028 (2011).

\bibitem{GL99}
A.M. Grundland and D. Levi, {\it On higher-order Riccati equations as B\"acklund transformations}, J. Phys. A
{\bf 32}, 3931--3937 (1999).

\bibitem{Palais}
R.S. Palais, {\sl Global Formulation of the Lie Theory of Transformation Groups}, Mem. Amer. Math. Soc. {\bf 22},
1957.

\bibitem{Ev90}
N.W. Evans, {\it Super-integrability of the Winternitz system}, Phys. Lett. A {\bf 147}, 483--486  (1990).

\bibitem{Ev91}
N.W. Evans, {\it Group theory of the Smorodinsky-Winternitz system}, J. Math. Phys. {\bf 32}, 3369--3375
(1991).

\bibitem{FMSU65}
J. Fri\v{s}, V. Mandrosov, Y.A. Smorodinsky, M. Uhl\'i\v{r} and P. Winternitz,
{\it On higher symmetries in quantum
mechanics},
Phys. Lett. {\bf 16}, 354--356 (1965).

\bibitem{GPS95}
C. Grosche,\! G.S.\! Pogosyan and A.N. Sissakian, {\it \!Path\! integral discussion for
Smorodinsky--Winternitz potentials I, two-dimensional and three-dimensional Euclidean space},
Fortschr. Phys. {\bf 43}, 453--521 (1995).

\bibitem{Pi50}
E. Pinney, {\it The nonlinear differential equation $\ddot y+p(x)y+cy^{−3} =0$},
Proc. A.M.S. {\bf 1}, 681 (1950).

\bibitem{Ka98}
R.S. Kaushal,
{\it Construction of exact invariants for time
dependent classical dynamical systems},
Int. J. Theor. Phys. {\bf 37}, 1793--1856 (1998).

\bibitem{DW10}
J. D'Ambroise and F.L. Williams,
{\it A dynamic correspondence between Bose–Einstein condensates and Friedmann–Lema\^itre–Robertson–Walker and Bianchi I cosmology with a cosmological constant},
J. Math. Phys. {\bf 51}, 062501 (2010).

\bibitem{LA08}
P.G.L. Leach, K. Andriopoulos,
{\it The Ermakov equation: a commentary},
Appl. Anal. Discrete Math. {\bf 2}, 146--157 (2008).

\bibitem{Re99}
R. Redheffer,
{\it Steen's equation and its generalizations},
Aequationes Math. {\bf 58}, 60--72 (1999).

\bibitem{BR97}
L.M. Berkovich and N.H. Rozov,
{\it Transformations of linear differential equations of second order and adjoined nonlinear equations},
Arch. Math. (Brno) {\bf 33} 75--98 (1997).

\bibitem{CFLSV}
J.F. Cari\~nena, J. de Lucas and C. Sard\'on,
{\it Lie-Hamilton systems:
theory and applications}, preprint 2012.

\bibitem{BM09II}
D. Bl\'azquez-Sanz and J.J. Morales-Ruiz,
{\it Local and global aspects of Lie superposition theorem}, J.
Lie Theory {\bf 20}, 483--517 (2010).

\bibitem{CRS05}
J.F. Cari{\~n}ena, M.F. Ra{\~n}ada and M. Santander,
{\it Lagrangian formalism for nonlinear second-order
Riccati systems: one-dimensional integrability and two-dimensional superintegrability}, J. Math. Phys.
\textbf{46}, 062703 (2005).

\bibitem{KL09}
A. Karasu and P.G.L. Leach,
{\it Nonlocal symmetries and integrable ordinary differential equations:
  $\ddot x + 3x\dot x + x^3 = 0$ and its generalizations},
J. Math. Phys. \textbf{50}, 073509 (2009).

\bibitem{Lie1880}
S. Lie,
{\sl Sophus Lie's 1880 transformation group paper},
Math. Sci. Press., Brookline, 1975.

\bibitem{GLKP90}
A. Gonz\'alez-L\'opez, N. Kamran and P. Olver,
{\it Lie algebras of vector fields in the real plane},
Proc. London Math. Soc. {\bf 64}, 339--368 (1992).

\bibitem{Da97}
R.W.R. Darling,
{\it Converting matrix Riccati equations to second-order linear ODE},
SIAM Rev. {\bf 39}, 508--510 (1997).


\end{thebibliography}
\end{document}